\documentclass[a4paper]{amsart}  

\usepackage{amssymb} \usepackage{amscd} \usepackage{times}

\def\zbb{\mathbb{Z}}  
  
  \def\phi{\varphi}
 \def\p1{{\mathbb{P}^1_\zbb}}

\newcommand{\be} {\begin{equation}}
\newcommand{\ee} {\end{equation}}

\begin{document}

\title{An existence result}
\author{Samy Skander Bahoura}
\address{Equipe d'Analyse Complexe et G\'eom\'etrie, Universit\'e Pierre et marie Curie, 75005, Paris, France}

\maketitle
\begin{abstract}
On compact Riemannian manifold of dimension $ n $, and under some conditions on the curvature, we have a changing-sign solutions for $ n $ large enough.\end{abstract}
\underbar{\bf Introduction}

\bigskip

Let $ (M, g) $ be a compact Riemannian manifold without boundary of
dimension $ n \geq 3 $. We consider the following equation:

$$ \Delta u = |u|^{4/(n-2)} u, \,\,\, u \not \equiv 0 \qquad (E) $$

Where $ \Delta = -\nabla^i(\nabla_i) $ is the Laplace-Beltrami operator on
$ M $.

\bigskip

{\it Holcman's Problem:} Is there a changing-sign solution to the problem $
(E) $ ?

\bigskip

In his paper, see [4], Holcman proved that, if the scalar
curvature $ R $ of $ M $ is positive somewhere, ($ R(P) > 0, P \in
M $), then $ (E) $ has a changing-sign solution.

\bigskip

Here, we have,

\bigskip

{\bf Theorem.} {\it Assume that $ (M, g) $ is not conformally flat
manifold of dimension $ n \geq 13 $ and $ R \equiv 0 $, then, $ (E) $
has a changing-sign solution}.

\bigskip

For the proof of the Theorem, we use T. Aubin's and Holcman's
methods and ideas and their computations, see [3,4]. We use the variational method with an explicit expansion of the Yamabe type functional. The Aubin's and Holcman's approach is the subcritical approach, they solve the nodal problem with subcritical exponent and they prove that the sequence with subcritical exponent converge to a solution of the problem (non-concentration), and the goal is to find tests functions which satisfy the inequality of non-concentration, linked to the Sobolev embedding, see also [2].

\bigskip

\underbar {\it Question:} Is it possible, if we use Schoen
inequality, in the proof of the Yamabe problem for the conformally
flat case, to have the same result for $ R \equiv 0 $ ?

\bigskip

\underbar{ \it Remark:} Look also, the paper of E. Humbert and B. Ammann, [1],
about " The second Yamabe invariant".

In our result, we can not do a conformal change of metric because we will have a term which contain the scalar curvature and in this case we change the equation $ (E) $.

\bigskip

\underbar{\bf Proof of the Theorem.}

\bigskip

Let us consider $ (M,g) $ a compact Riemannian manifold without
boundary and not locally conformally flat. Assume that the scalar
curvature $ R \equiv 0 $ and we work with the Yamabe functional:

$$ J(\phi)=\dfrac{\int_M |\nabla(\phi)|^2 dV_g+\int_M R\phi^2
dV_g}{\left (\int_M \phi^N \right )^{2/N}}=\dfrac{\int_M
|\nabla(\phi)|^2 dV_g}{\left (\int_M \phi^N \right )^{2/N}}. $$

Let $ P $ the point where $ Weyl_g(P)\not = 0 $,

As in the paper of T. Aubin, we do a conformal change of metric $
\tilde g=\psi^{4/(n-2)} g $ such that:

$$ \tilde J(\phi_{\epsilon})=\dfrac{1}{K}[1-|Weyl_{\tilde g}(P)|^2\epsilon^4+o(\epsilon^4)], $$

where $ \tilde J $ is the Yamabe functional for the metric $ \tilde
g $ and $ \phi_{\epsilon} $ the following functions:

$$ \phi_{\epsilon}(\tilde r)=\dfrac{\epsilon^{(n-2)/2}}{(\epsilon^2+{\tilde r}^2)^{(n-2)/2}}-\dfrac{\epsilon^{(n-2)/2}}{(\epsilon^2+{\tilde \delta }^2)^{(n-2)/2}} \,\,\, {\rm if} \,\, \tilde r=\tilde d(P,x) \leq \tilde \delta,\,\,\, {\rm otherwise } \,\, 0. $$

Also, we know that:

$$ J(\psi \phi_{\epsilon})=\tilde J(\phi_{\epsilon}) . $$

Let us consider the following functions:

$$ \bar \phi_{\epsilon}=\psi (\phi_{\epsilon}-\mu_{\epsilon}), $$

with, $ \mu_{\epsilon} >0 $ is such that:

$$ \int_M |\psi(\phi_{\epsilon}-\mu_{\epsilon})|^{N-2}[\psi(\phi_{\epsilon}-\mu_{\epsilon})]dV_g=0. $$

If we compute with $ \tilde g $, we have:

$$ \int_M \dfrac{1}{\psi}|\phi_{\epsilon}-\mu_{\epsilon}|^{N-2}(\phi_{\epsilon}-\mu_{\epsilon})d \tilde V=\int_M f|\phi_{\epsilon}-\mu_{\epsilon}|^{N-2}(\phi_{\epsilon}-\mu_{\epsilon})d \tilde V=0. $$

with $ f=\dfrac{1}{\psi} >0 $.

We know, see Holcman, that $ \mu_{\epsilon} $ is equivalent to $
{\epsilon}^{[(n-2)^2]/2(n+2)} $ for $ \epsilon $ near $ 0 $.

\bigskip

Since the distance function $ \tilde r $ is Lipschitzian and equivalent to the first distance function $ r $, we can compute (when we have the gradient), with respect to the $ \tilde r $. We can write,

$$ \int_M |\nabla [\psi(\phi_{\epsilon}-\mu_{\epsilon})]|^2 \leq \int_M |\nabla (\psi \phi_{\epsilon})|^2dV_g + c_1\mu_{\epsilon}, $$

to see this, we write:

$$  \int_M |\nabla [\psi(\phi_{\epsilon}-\mu_{\epsilon})]|^2 = \int_M |\nabla (\psi \phi_{\epsilon})|^2dV_g + 2\mu_{\epsilon} \int_M <\nabla \psi, \nabla (\psi \phi_{\epsilon})> + O(\mu_{\epsilon}^2), $$

$$ \int_M <\nabla \psi, \nabla (\psi \phi_{\epsilon})> = O(\int_M \phi_{\epsilon})+ O(\int_M | \nabla (\phi_{\epsilon})|), $$  

We can see that;

$$ \tilde \nabla^i(\phi_{\epsilon})=\psi^{-4/(n-2)} \nabla^i (\phi_{\epsilon}), $$

Thus, for two positive constants $ C_1, C_2 $, we have:

$$ C_2|\tilde \nabla \phi_{\epsilon}|\leq |\nabla \phi_{\epsilon}|\leq C_1|\tilde \nabla \phi_{\epsilon}|=C_1|\partial_{\tilde r}\phi_{\epsilon}(\tilde r)|, $$


$$ \int_M | \nabla \phi_{\epsilon}| dV_g  \leq C_4 \int_M|\partial_{\tilde r}\phi_{\epsilon}(\tilde r)| d \tilde V \leq C_5, $$

and,

$$ \int_M | \nabla (\psi\phi_{\epsilon})|^2 dV_g= \int_M | \tilde \nabla \phi_{\epsilon}|^2 d\tilde V + o(1)\geq C_6 \int_M|\partial_{\tilde r}\phi_{\epsilon}(\tilde r)|^2 d \tilde V \geq C_7 >0, $$

(see, Aubin computations), and we have the result for the gradient.

And, we have:

$$ (\int_M |\psi(\phi_{\epsilon}-\mu_{\epsilon})|^N)^{2/N}\geq (\int_M |\psi \phi_{\epsilon}|^N)^{2/N}-c_2\mu_{\epsilon} , \,\, {\rm with } \,\, c_2>0. $$

because,

$$ ||\psi (\phi_{\epsilon}-\mu_{\epsilon})||_{L^N, g}^N = \int_M |\psi(\phi_{\epsilon}-\mu_{\epsilon})|^N dV=\int_M |(\phi_{\epsilon}-\mu_{\epsilon})|^N d \tilde V= ||(\phi_{\epsilon}-\mu_{\epsilon})||_{L^N,  \tilde g}^N, $$

and, for $ \tilde g $

$$ ||\phi_{\epsilon}||_{L^N, \tilde g} \leq ||(\phi_{\epsilon}-\mu_{\epsilon})||_{L^N,  \tilde g} + |M|^{1/N}\mu_{\epsilon}, $$

and, because $ ||(\phi_{\epsilon}-\mu_{\epsilon})||_{L^N, \tilde g} \to c >0 $ (or,  $ ||\phi_{\epsilon}||_{L^N, \tilde g} \to c' >0 $,) (see the computations of Holcman's paper with the metric $ \tilde g $),

$$ ||\psi\phi_{\epsilon}||_{L^N, g}^2= ||\phi_{\epsilon}||_{L^N, \tilde g}^2 \leq ||(\phi_{\epsilon}-\mu_{\epsilon})||_{L^N, \tilde g}^2 + c_2\mu_{\epsilon}, $$

and then,

$$ \dfrac{\int_M |\nabla [\psi(\phi_{\epsilon}-\mu_{\epsilon})]|^2}{\left ( \int_M |\psi(\phi_{\epsilon}-\mu_{\epsilon})|^N \right )^{2/N}}\leq J(\psi \phi_{\epsilon})(1+c_3\mu_{\epsilon}) $$

Thus,

$$ \dfrac{\int_M |\nabla [\psi(\phi_{\epsilon}-\mu_{\epsilon})]|^2}{\left ( \int_M |\psi(\phi_{\epsilon}-\mu_{\epsilon})|^N \right )^{2/N}} \leq \dfrac{1}{K}[1+c_4\epsilon^{(n-2)^2/2(n+2)}-|Weyl_{\tilde g}(P)|^2 \epsilon^4 + o(\epsilon^4) ]. $$

We can say that,  $ \epsilon^{(n-2)^2/2(n+2)} $ is very small if we
compare it to $ \epsilon^4 $ if, $ \dfrac{(n-2)^2}{2(n+2)} > 4 $, and
then, if $ n\geq 13 $.

Thus, on $ M $, we have test functions 

$$ \bar \phi_{\epsilon}=\psi (\phi_{\epsilon}-\mu_{\epsilon}) \not \equiv 0, $$

such that:

$$ \int_M |\bar \phi_{\epsilon}|^{N-2} \bar \phi_{\epsilon} dV_g=0, $$

and, the Sobolev quotient is such that:

$$ \dfrac{\int_M |\nabla \bar \phi_{\epsilon}|^2}{\left ( \int_M |\bar \phi_{\epsilon}|^N \right )^{2/N}} \leq \dfrac{1}{K}[1+c_4\epsilon^{(n-2)^2/2(n+2)}-|Weyl_{\tilde g}(P)|^2 \epsilon^4 + o(\epsilon^4) ] < \dfrac{1}{K}. $$
 
Thus, the variational problem has a nodal solution on $ M $.

\bigskip

\underbar{ \it Remark 1:} We can replace $ \mu_{\epsilon} $ by $\mu_{\epsilon}^2 $, in this case we can assume $ n \geq 9 $.

\smallskip

\underbar{ \it Remark 2:} This method works if we assume that, there is a point  $ P $ such that $ Weyl_g (P) \not = 0 $ and $ R \equiv 0 $ in the neighborhood of $ P $. (Such manifolds exist, it is sufficient to solve the prescribed scalar curvature problem for non-positive scalar curvature, by considering the condition on the first eigenvalue of small balls, see Rauzy and Veron in  [1]).

\smallskip

\underbar{ \it Remark 3:} This method works if we assume that, there is a point  $ P $ such that $ Weyl_g (P) \not = 0 $ and $ R(P)=\nabla R(P)= \nabla^2R(P)=0 $.

\smallskip

We have the following corollary:

\smallskip

{\bf Corollary.} {\it Assume that $ (M, g) $ is not conformally flat manifold of dimension $ n \geq 13 $ and $ R \equiv 0 $ in a neighborhood of a point $ P $ such that $ Weyl_g(P) \not = 0 $, then, $ (E) $
has a changing-sign solution}.

\end{document}